\documentclass{amsart}
\usepackage{mathptmx, amscd, amssymb, enumerate, colonequals, stmaryrd, xcolor}
\usepackage{color}
\definecolor{chianti}{rgb}{0.6,0,0}
\definecolor{meretale}{rgb}{0,0,.6}
\definecolor{leaf}{rgb}{0,.35,0}
\usepackage[colorlinks=true, pagebackref, hyperindex, citecolor=meretale, urlcolor=leaf, linkcolor=chianti]{hyperref}

\newtheorem{theorem}{Theorem}[section]

\theoremstyle{definition}

\newtheorem{remark}[theorem]{Remark}
\newtheorem{question}[theorem]{Question}
\numberwithin{equation}{theorem}

\def\depth{\operatorname{depth}}
\def\ffrac{\operatorname{frac}}
\def\tr{{\operatorname{tr}}}
\def\Hom{\operatorname{Hom}}
\def\GL{\operatorname{GL}}

\def\frakA{\mathfrak{A}}
\def\frakC{\mathfrak{C}}
\def\frakP{\mathfrak{P}}

\def\FF{\mathbb{F}}
\def\FFp{\mathbb{F}_{\!p}}
\def\NN{\mathbb{N}}
\def\ZZ{\mathbb{Z}}

\def\ph{\phantom{p}}
\def\ge{\geqslant}
\def\le{\leqslant}
\def\to{\longrightarrow}
\def\mapsto{\longmapsto}
\def\into{\lhook\joinrel\longrightarrow}

\begin{document}
\title{Abelian extensions of equicharacteristic regular rings need not be Cohen-Macaulay}

\author{Aryaman Maithani}
\address{Department of Mathematics, University of Utah, Salt Lake City, Utah 84112}
\email{maithani@math.utah.edu}

\author{Anurag K. Singh}
\address{Department of Mathematics, University of Utah, Salt Lake City, Utah 84112}
\email{singh@math.utah.edu}

\author{Prashanth Sridhar}
\address{Department of Mathematics, University of Alabama, Tuscaloosa, Alabama 35487}
\email{psridhar1@ua.edu}

\thanks{A.M. was supported by a Simons Dissertation Fellowship; A.M. and A.K.S. were supported by NSF grants DMS~2101671 and DMS~2349623.}

\dedicatory{To Professor Paul Roberts, on the occasion of his eightieth birthday}

\begin{abstract}
By a theorem of Roberts, the integral closure of a regular local ring in a finite abelian extension of its fraction field is Cohen-Macaulay, provided that the degree of the extension is coprime to the characteristic of the residue field. We show that the result need not hold in the absence of this requirement on the characteristic: for each positive prime integer $p$, we construct polynomial rings over fields of characteristic $p$, whose integral closure in an elementary abelian extension of order $p^2$ is not Cohen-Macaulay. Localizing at the homogeneous maximal ideal preserves the essential features of the construction.
\end{abstract}
\maketitle

%%%%%%%%%%%%%%%%%%%%%%%%%%%%%%%%%%%%%%%%%%%%%%%%%%%%%%%%%%%%%%%%%%%%%%%%
\section{Introduction}
%%%%%%%%%%%%%%%%%%%%%%%%%%%%%%%%%%%%%%%%%%%%%%%%%%%%%%%%%%%%%%%%%%%%%%%%

Paul Roberts \cite{Roberts} proved the following:

\begin{theorem}[Roberts]
\label{theorem:roberts}
Let $R$ be a regular local ring with fraction field $K$. Let $L$ be a finite Galois extension of $K$ with an abelian Galois group. Assume moreover that the order of the Galois group is not divisible by the characteristic of the residue field of $R$. Let $S$ denote the integral closure of $R$ in $L$. Then~$S$ is a Cohen-Macaulay ring.
\end{theorem}

If one instead assumes that $R$ is a UFD rather than a regular local ring, with the other hypotheses still in place, Roberts observes that the proof yields that $S$ is a free $R$-module; if the requirement on the characteristic is dropped, \cite[Example~1]{Roberts} provides a UFD $R$ of mixed characteristic $2$, and an extension field $L$ of degree $2$ over $\ffrac(R)$, such that if $S$ denotes the integral closure of $R$ in $L$, then $S$ is not a free $R$-module; in fact, no nonzero $S$-module is free over $R$.

Returning to the case where $R$ is regular, the abelian hypothesis cannot be weakened to solvable or nilpotent in view of \cite[Example~2]{Roberts}. Additionally, Koh demonstrated that Theorem~\ref{theorem:roberts} may fail in the absence of the requirement on the order of the group: \cite[Example~2.4]{Koh} is an example of a regular local ring $R$ of mixed characteristic $3$, with an extension field~$L$ having Galois group $\ZZ/3\ZZ$ over $\ffrac(R)$, such that the integral closure of~$R$ in $L$ is not Cohen-Macaulay. Nevertheless, there has been work exploring extensions of Roberts's theorem by two distinct methods: in terms of tracking ramification in codimension one, and in terms of constructing birational maximal Cohen-Macaulay modules, see \cite{Katz, Katz:Sridhar, Sridhar1, Sridhar2}. The purpose of this brief note is to record equal characteristic examples where the integral closure of a regular ring in an elementary abelian extension of order $p^2$ is not Cohen-Macaulay---indeed, the Cohen-Macaulay defect can be arbitrarily large---and to demonstrate that such examples arise quite naturally in the context of modular invariant theory. While the examples are recorded in the framework of $\NN$-graded rings, the relevant properties are preserved upon localization at the respective homogeneous maximal ideals.

%%%%%%%%%%%%%%%%%%%%%%%%%%%%%%%%%%%%%%%%%%%%%%%%%%%%%%%%%%%%%%%%%%%%%%%%
\section{Modular invariant rings}
%%%%%%%%%%%%%%%%%%%%%%%%%%%%%%%%%%%%%%%%%%%%%%%%%%%%%%%%%%%%%%%%%%%%%%%%

Let $V$ be a finite rank vector space over a field; an element of $\GL(V)$ is a \emph{bireflection} if it fixes a subspace of $V$ of codimension at most $2$. The following is \cite[Corollary~3.7]{Kemper:bireflection}:

\begin{theorem}[Kemper]
\label{theorem:kemper}
Let $K$ be a field of characteristic $p>0$, and $H$ a finite subgroup of~$\GL_n(K)$, with its natural action on the polynomial ring $T\colonequals K[x_1,\dots,x_n]$.

If $H$ is a $p$-group, and $T^H$ is Cohen-Macaulay, then $H$ is generated by bireflections.
\end{theorem}

Using this, we prove:

\begin{theorem}
\label{theorem:main}
For each positive prime integer $p$, there exists a polynomial ring~$R$ over the finite field~$\FFp$, and an elementary abelian extension $L$ of $\ffrac(R)$ of order $p^2$, such that the integral closure of $R$ in $L$ is not Cohen-Macaulay.
\end{theorem}

\begin{proof}
The ring $R$ will be constructed as the invariant ring $T^G$ for the action of a group $G$ on a polynomial ring $T$, while the extension ring $S$ that is not Cohen-Macaulay equals $T^H$, for a subgroup $H$ of $G$. Schematically,
\[
\CD
T\colonequals \FFp[x_1,y_1,x_2,y_2,x_3,y_3] \\
\lvert \\
S\colonequals T^H \\
\lvert \\
R\colonequals T^G,
\endCD
\]
with $R$ being a polynomial ring.

Let $G$ be the subgroup of $\GL_6(\FFp)$ consisting of the matrices
\[
\phantom{\quad\text{ for } a,b,c\in\FFp}
\left[
\begin{array}{cc|cc|cc}
1 & a & & & &\\
& 1 & & & &\\
\hline
& & 1 & b & &\\
& & & 1 & &\\
\hline
& & & & 1 & c\\
& & & & & 1
\end{array}
\right]
\quad\text{ for } a,b,c\in\FFp.
\]
The group $G$ is abelian, isomorphic to the direct product of three copies of $\ZZ/p\ZZ$. Consider the linear action of $G$ on the polynomial ring $T$, where $M\in G$ acts via
\[
M\colon X\mapsto MX,
\]
with $X\colonequals(x_1,y_1,x_2,y_2,x_3,y_3)^\tr$ denoting the column vector of indeterminates. We claim that~$T^G$, i.e., the invariant ring for the action of $G$ on $T$, is the polynomial ring
\[
R\colonequals\FFp[y_1^{\ph},\ x_1^p-x_1^{\ph}y_1^{p-1},\ y_2^{\ph},\ x_2^p-x_2^{\ph}y_2^{p-1},\ y_3^{\ph},\ x_3^p-x_3^{\ph}y_3^{p-1}].
\]
It is readily seen that $R\subseteq T^G$. For each $i$, the element $x_i\in T$ is a root of the polynomial
\[
Z^p - Zy_i^{p-1} - \big(x_i^p-x_i^{\ph} y_i^{p-1}\big)\ \in\ R[Z],
\]
so $T$ is integral over $R$ and $[\ffrac(T):\ffrac(R)] \le p^3$. Since $G$ has order $p^3$, it follows that
\[
\ffrac(R)\ =\ \ffrac(T)^G \ =\ \ffrac(T^G),
\]
where the second equality holds since $G$ is finite. Since $R$ has dimension $6$ and is generated by the $6$ elements displayed, $R$ is a polynomial ring, hence normal. Each element of~$T^G$ lies in $\ffrac(R)$ and is integral over $R$, so the normality of $R$ yields that $R=T^G$, as claimed.

We note that while one direction of the Shephard-Todd-Chevalley-Serre theorem may fail in the modular case, \cite[page~3]{Serre}, one may instead use~\cite[Theorem~1.4]{Nakajima} or Remark~\ref{remark:general} below to conclude that the invariant ring $T^G$ is a polynomial ring; this is in lieu of the direct proof that we have chosen to include above.

Next, consider the cyclic subgroup $H$ of $G$ generated by the matrix
\[
\left[
\begin{array}{cc|cc|cc}
1 & 1 & & & &\\
& 1 & & & &\\
\hline
& & 1 & 1 & &\\
& & & 1 & &\\
\hline
& & & & 1 & 1\\
& & & & & 1
\end{array}
\right].
\]
Since $H$ is a $p$-group that is \emph{not} generated by bireflections---the only bireflection in $H$ is the identity---Theorem~\ref{theorem:kemper} implies that the invariant ring $S\colonequals T^H$ is not Cohen-Macaulay.

To summarize, we have constructed a polynomial ring $R$ with a finite normal extension ring $S$ that is not Cohen-Macaulay; the Galois group of $\ffrac(S)$ over $\ffrac(R)$ is
\[
G/H\ \cong \ \ZZ/p\ZZ\times\ZZ/p\ZZ,
\]
the elementary abelian group of order $p^2$, and $S$ is the integral closure of $R$ in $\ffrac(S)$.
\end{proof}

\begin{remark}
The ring $S$ in the proof of Theorem~\ref{theorem:main} is a UFD since there is no nontrivial character $H\to\FFp^\times$, see, for example, \cite[Corollary~3.9.3]{Benson}. While $S$ is not Cohen-Macaulay, it evidently admits a small Cohen-Macaulay algebra, namely the polynomial ring $T$. Consequently, $T$ is not a free $S$-module. 
\end{remark}

\begin{remark}
\label{remark:general}
Extending the construction in the proof of Theorem~\ref{theorem:main}, let $P$ be a nontrivial abelian $p$-group acting faithfully on a finite rank $\FFp$-vector space $V$, such that the invariant ring $\FFp[V]^P$ is polynomial; this holds, for example, if the fixed subspace $V^P$ has codimension one, see \cite[Theorem~3.9.2]{Campbell:Wehlau}. For an integer $d\ge 3$, consider the~$d$-fold product $G\colonequals P\times\dots\times P$ acting on the polynomial ring $T\colonequals \FFp[V^{\oplus d}]$. Then the invariant ring $T^G$ is the tensor product of $d$ copies of~$\FFp[V]^P$, hence a polynomial ring. Set $H$ to be the diagonal copy of $P$ in $G$, which contains no bireflections other than the identity, since $d\ge 3$. One obtains a tower of rings
\[
\CD
\FFp[V^{\oplus d}] \\
\lvert \\
S\colonequals {\FFp[V^{\oplus d}]}^H \\
\lvert \\
R\colonequals {\FFp[V^{\oplus d}]}^G,
\endCD
\]
where $S$ is not Cohen-Macaulay in view of Kemper's result, Theorem~\ref{theorem:kemper}. The Galois group of $\ffrac(S)$ over $\ffrac(R)$ is the abelian group $G/H$.

Specifically, suppose $P$ is the cyclic group generated by
\[
\begin{bmatrix}
1 & 1 \\
0 & 1 
\end{bmatrix},
\]
with its usual action on the rank two vector space $\FFp^2$. Let $G$ be the product of $d$ copies of~$P$ for an integer~$d\ge 3$, and let $T$ and $H$ be as above; note that $\dim T=2d$. The fact that the invariant ring $S\colonequals T^H$ is not Cohen-Macaulay also follows from the work of Ellingsrud and Skjelbred: \cite[Corollaire~2.4]{ES} implies that $S$ does not satisfy the Serre condition $(S_3)$, while \cite[Corollaire~3.2]{ES} implies that $S$ has depth $d+2$. In particular, the \emph{Cohen-Macaulay defect} of the ring $S$, i.e., $\dim S-\depth S$, is $2d-(d+2)=d-2$.
\end{remark}

While we have constructed examples with extension degree $p^2$, the following remains:

\begin{question}
\label{question:cyclic}
Let $R$ be a regular local ring of characteristic $p>0$, with fraction field $K$. Let $L$ be a Galois extension of $K$ with Galois group $\ZZ/p\ZZ$, and $S$ the integral closure of $R$ in $L$. Is $S$ a Cohen-Macaulay ring?
\end{question}

The analogous statement is false in mixed characteristic by~\cite[Example~2.4]{Koh}. In view of Roberts' result, as discussed in the introduction, one may ask: Let $R$ be a local \emph{UFD} of characteristic $p>0$, with fraction field $K$. Let $L$ be a Galois extension of $K$ with Galois group $\ZZ/p\ZZ$, and $S$ the integral closure of $R$ in $L$. Is $S$ a \emph{free} $R$-module? The answer to this is negative: the ring $T^H$ from Theorem~\ref{theorem:main} is a UFD, while $T$ is not a free $T^H$-module. 

We remark that the answer to Question~\ref{question:cyclic} is positive for $p=2$: the inclusion $R \into S$ splits as $R$ is regular, so $S$ is isomorphic to $R \oplus M$ for $M$ a reflexive rank one $R$-module; but then $M$ is isomorphic to $R$ since $R$ is a UFD\@. We can also answer Question~\ref{question:cyclic} in the affirmative in the following special case:

\begin{remark}
Continuing with the notation as in Question~\ref{question:cyclic}, the field extension $L/K$ is given by an Artin-Schreier polynomial $X^p-X-a/b \in K[X]$, with $a,b\in R$ coprime. We show that $S$ is Cohen-Macaulay if $b$ is a squarefree element of $R$. We may assume $p\ge 3$.

Clearing denominators, $L/K$ is also the field extension corresponding to
\[
F(X) \colonequals X^p-b^{p-1}X-ab^{p-1}\ \in\ R[X],
\]
and $S$ is the integral closure of $A\colonequals R[X]/(F)$. Recall that if $\frakA \subseteq \ffrac(A)$ is a fractional ideal, its inverse $\frakA^{-1}$ is $\{f \in \ffrac(A) \mid f\frakA \subset A\}$, and this is isomorphic to $\Hom_A(\frakA, A)$. Consider the conductor $\frakC \colonequals S^{-1}$ of the extension $A \subseteq S$. If $\frakC$ is the unit ideal, we are done; assume that is not the case. Since $A$ is Gorenstein, \cite[Proposition~2.11]{Sridhar2} implies that~$S$ is Cohen-Macaulay if and only if $A/\frakC$ is so. Towards this, we next compute the primary decomposition of $\frakC$ in $A$.

Let $\omega$ denote the image of $X$ in $A$, so $F(\omega) = 0$. Let $\frakP$ be an associated prime of~$\frakC$ in~$A$. Since~$A$ is $(S_2)$ and $\frakC$ is an $A$-dual, $\frakC$ is $(S_2)$, so~$\frakP$ has height one. The separability of $L/K$ implies that $F'(\omega) = -b^{p-1} \in \frakC$, see, for example, \cite[Theorem~12.1.1]{SH}. Thus,~$\frakP$ contains a prime factor $\pi$ of $b$ in $R$. Moreover,~$\frakP$ contains $\omega$ as $\omega^p = b^{p-1}(\omega + a)$. But~$(\omega, \pi)A$ is a prime ideal of height one, thus it equals $\frakP$. Hence each associated prime of~$\frakC$ is of the form $\frakP = (\omega,\pi)A$, for some prime factor $\pi$ of $b$ in $R$. We next compute $\frakC_\frakP$.

As $b$ is squarefree, we may rewrite $F = X^p - b^{p-1}(X + a)$ as $F = X^p - \pi^{p-1} \lambda$ for some unit $\lambda \in R[X]_{(X, \pi)}$. 
The element $\pi/\omega$ is integral over $A_\frakP$ since $(\pi/\omega)^{p-1} = \lambda^{-1}\omega$, so
\[
A_\frakP[\pi/\omega]\ \subseteq\ S_\frakP.
\]
We claim that equality holds above. Each maximal ideal of $A_\frakP[\pi/\omega]$ contains $\omega$, $\pi$, and also $\pi/\omega$, using once again that $(\pi/\omega)^{p-1} = \lambda^{-1}\omega$. But $(\pi, \omega, \pi/\omega)A_\frakP[\pi/\omega]$ is a maximal ideal, so $A_\frakP[\pi/\omega]$ is local, moreover with maximal ideal generated by $\pi/\omega$. It follows that $A_\frakP[\pi/\omega]$ is a regular local ring, hence normal, proving the claim. Note that $S_\frakP$ is a discrete valuation ring with uniformizer $\pi/\omega$, while the maximal ideal $(\omega, \pi)A_\frakP$ of the subring $A_\frakP$ is generated by
\[
(\pi/\omega)^{p-1}=\lambda^{-1}\omega\quad\text{ and }\quad (\pi/\omega)^p=\lambda^{-1}\pi.
\]
Upon completion, the inclusion $A_\frakP \subset S_\frakP$ takes the form
\[
\FF\llbracket(\pi/\omega)^{p-1},\, (\pi/\omega)^p\rrbracket\ \subset\ \FF\llbracket\pi/\omega\rrbracket,
\]
for $\FF$ a coefficient field, so the calculation of the conductor is essentially the calculation of the Frobenius number of the semigroup generated by $p-1$ and $p$, and this is $p^2-3p+1$. It follows that $\frakC_\frakP = \frakP^{p-2} A_\frakP$, and, in turn, that $\frakC A_\frakP \cap A$ is the symbolic power $\frakP^{(p-2)}$.

As $\frakC$ is unmixed, we obtain the primary decomposition $\frakC = \bigcap_i (\omega,\pi_i)^{(p-2)}$, where $\pi_i$ ranges over the prime factors of $b$ in $R$. Note that $F \in (X, \pi_i)^{p-2}$ for each $i$, and thus,
\[
A/\frakC\ \cong\ \frac{R[X]}{\bigcap_i (X,\pi_i)^{(p-2)}}\ =\ \frac{R[X]}{\bigcap_i (X,\pi_i)^{p-2}},
\] 
where the equality holds as ordinary and symbolic powers coincide for ideals generated by a regular sequence. Moreover, in the regular ring $R[X]$, it is readily seen that
\[
\bigcap\nolimits_i (X,\pi_i)^{p-2}\ =\ (X,b)^{p-2}
\]
since $b$ is squarefree, with the elements $\pi_i$ as its prime factors. Since $(X,b)^{p-2}$ is a perfect ideal of $R[X]$, we conclude that the ring $A/\frakC \cong R[X]/(X,b)^{p-2}$ is Cohen-Macaulay.
\end{remark}

%%%%%%%%%%%%%%%%%%%%%%%%%%%%%%%%%%%%%%%%%%%%%%%%%%%%%%%%%%%%%%%%%%%%%%%%

\end{document}